\documentclass{amsart}
\usepackage{amssymb}
\usepackage{eucal}

\theoremstyle{plain}
\newtheorem{theorem}{Theorem}
\newtheorem{THEOREM}{Theorem}

\newtheorem{proposition}[theorem]{Proposition}

\newtheorem{corollary}[theorem]{Corollary}

\theoremstyle{definition}
\newtheorem{definition}[theorem]{Definition}
\newtheorem*{definition*}{Definition}

\newtheorem*{notation*}{Notation}

\theoremstyle{remark}

\newtheorem*{remark*}{Remark}

\newtheorem*{remarks*}{Remarks}

\newtheorem*{example*}{Example}

\newtheorem*{examples*}{Examples}

\renewcommand\theenumi{{\roman{enumi}}}

\newcommand\iso{{\widetilde\to}}                    %Invertible morphism

\newcommand\moprm[2]{\newcommand{#1}{\mathop{\mathrm{#2}}\nolimits}}     %New math operator in roman font
\newcommand\mopit[2]{\newcommand{#1}{\mathop{\mathit{#2}}\nolimits}}     %New math operator in italic font
\newcommand\mopbf[2]{\newcommand{#1}{\mathop{\mathbf{#2}}\nolimits}}     %New math operator in bold font

\moprm{\Hom}{Hom}                                   %Homomorphisms (set)
\moprm{\End}{End}                                   %Endomorphisms (set)
\moprm{\Aut}{Aut}                                   %Automorphisms (set)
\moprm{\Iso}{Iso}                                   %Isomorphisms (set)
\mopit{\HOM}{\cH om}                                %Homomorphisms (sheaf)
\mopit{\END}{\cE nd}                                %Endomorphisms (sheaf)
\mopit{\AUT}{\cA ut}                                %Automorphisms (sheaf)
\mopit{\ISO}{\cI so}                                %Isomorphisms (sheaf)

\moprm{\Ext}{Ext}                                   %Ext's (set)
\moprm{\Tor}{Tor}                                   %Tor's (set)
\mopit{\EXT}{\cE xt}                                %Ext's (sheaf)
\mopit{\TOR}{\cT or}                                %Tor's (sheaf)

\moprm{\ob}{Ob}                                     %The class of objects of a category

\moprm{\tr}{tr}                                     %Trace of an operator
\moprm{\rk}{rk}                                     %Rank of a bundle
\moprm{\ad}{ad}                                     %Adjoint action of a lie algebra
\moprm{\Ad}{Ad}                                     %Adjoint action of a lie group
\moprm{\id}{id}                                     %Identity morphism
\moprm{\supp}{supp}                                 %Support of a sheaf
\moprm{\chr}{char}                                  %Characteristic of a field
\moprm{\codim}{codim}                               %Codimension

\moprm{\spec}{Spec}                                 %Spectrum of an algebra
\moprm{\spf}{Spf}                                   %Formal spectrum of a topological algebra

\mopbf{\gm}{G_m}                                    %The multiplicative group
\mopbf{\ga}{G_a}                                    %The additive group
\moprm{\detrg}{detR\Gamma}
\moprm{\Pic}{Pic}

\newcommand{\oJ}{{\overline J}}
\newcommand{\oP}{{\overline P}}
\newcommand{\tC}{{\tilde C}}
\newcommand{\tJ}{{\tilde J}}
\newcommand{\tN}{{\tilde N}}
\newcommand{\cC}{{\mathcal C}}
\newcommand{\cF}{{\mathcal F}}
\newcommand{\cH}{{\mathcal H}}
\newcommand{\cJ}{{\mathcal J}}
\newcommand{\cK}{{\mathcal K}}
\newcommand{\cS}{{\mathcal S}}

\newcommand{\tg}{{\tilde g}}

\newcommand\kkk{{\Bbbk}}
\newcommand\cv{{\mathcal M}}
\newcommand\tpi{{\tilde \pi}}
\newcommand\frj{{\mathfrak{j}}}
\newcommand\Higgs{{\mathcal{H}iggs}}
\newcommand\SCurves{{\mathcal{SC}urves}}

\begin{document}

\title{Cohomology of line bundles on compactified Jacobians}
\author{D. Arinkin}

\begin{abstract} Let $C$ be an integral projective curve with planar singularities. For  the compactified
Jacobian $\oJ$ of $C$, we prove that topologically trivial line bundles on $\oJ$ are in one-to-one correspondence
with line bundles on $C$ (the autoduality conjecture), and compute the cohomology of $\oJ$ with coefficients in these
line bundles. We also show that the natural Fourier-Mukai functor from the derived category of quasi-coherent sheaves
on $J$ (where $J$ is the Jacobian of $X$) to that of quasi-coherent sheaves on $\oJ$ is fully faithful.
\end{abstract}

\maketitle

\section*{Introduction}

Let $C$ be a smooth irreducible projective curve over a field $\kkk$, and let $J$ be the Jacobian of $C$. As an abelian variety, $J$ is self-dual.
More precisely, $J\times J$ carries a natural line bundle (the Poincar\'e bundle) $P$
that is universal as a family of topologically trivial line bundles on $J$.

The Poincar\'e bundle defines the Fourier-Mukai functor
$${\mathfrak F}:D^b(J)\to D^b(J):\cF\mapsto Rp_{2,*}(p_1^*(\cF)\otimes P).$$
Here $D^b(J)$ is the derived category of quasi-coherent sheaves on $J$ and $p_{1,2}:J\times J\to J$ are the projections.
Mukai (\cite{Mukai}) proved that $\mathfrak F$ is an equivalence of categories; the proof uses the formula
\begin{equation}
Rp_{1,*}P\simeq O_\zeta[-g],
\label{eq:pushforward}
\end{equation}
where $O_\zeta$ is the structure sheaf of the zero element $\zeta\in J$ and $g$ is the genus of $C$.
Formula \eqref{eq:pushforward} goes back
to Mumford (see the proof of the theorem in \cite[Section III.13]{Mumford}).

Now suppose that $C$ is a singular curve, which we assume to be projective and integral.
The Jacobian $J$ is no longer projective, but it admits a natural compactification $\oJ\supset J$. By definition,
$\oJ$ is the moduli space of torsion-free sheaves $F$ on $C$ such that $F$
has generic rank one and $\chi(F)=\chi(O_C)$;
$J$ is identified with the open subset of locally free sheaves $F$.
It is natural to ask whether $\oJ$ is in some sense self-dual.
For instance, one can look for a Poincar\'e sheaf (or complex of sheaves)
$\overline P$ on $\oJ\times\oJ$. One can then ask whether $\overline P$ is,
in some sense, a universal family of sheaves
on $J$ and whether the corresponding Fourier-Mukai functor $\overline{\mathfrak F}:D^b(\oJ)\to D^b(\oJ)$ is an
equivalence.

In the case when singularities of $C$ are nodes or cusps, such Poincar\'e sheaf $\overline P$ is constructed
by E.~Esteves and S.~Kleiman in \cite{compactified}; they also prove the universality of $P$. In addition, if $C$
is a singular plane cubic, $\overline{\mathfrak F}$ is known to be an equivalence
(\cite{BurbanKreussler0,BurbanKreussler}, also formulated as Theorem 5.2 in \cite{BenZviNevins}).

If singularities of $C$ are more general, constructing the Poincar\'e sheaf $\oP$ on $\oJ\times\oJ$ is much harder (see 
Remark~(i) at the end of the introduction).
However, it is easy to construct a Poincar\'e bundle $P$ on $J\times\oJ$. It can then be used to define a Fourier-Mukai
transform
\begin{equation}
{\mathfrak F}:D^b(J)\to D^b(\oJ):\cF\mapsto Rp_{2,*}(p_1^*(\cF)\otimes P).
\label{eq:fouriermukai}
\end{equation}
Here it is important to work with the derived categories of quasicoherent sheaves, since
${\mathfrak F}$ does not preserve coherence. 

In this paper, we assume that $C$ is an
integral projective curve with planar singularities; the main result is that the formula \eqref{eq:pushforward}
still holds in this case. This implies that \eqref{eq:fouriermukai} is fully faithful. As a simple corollary, we
prove the following autoduality result: $P$ is the universal family
of topologically trivial  line bundles on $\oJ$, so that $J$ is identified with the connected component of the trivial
bundle in the moduli space of line bundles on $\oJ$. This generalizes the Autoduality Theorem of \cite{autoduality}
(see the remark after Theorem \ref{th:Jacobians}).

\begin{remarks*}
\setcounter{enumi}{1}
(\theenumi)\stepcounter{enumi}
Suppose that there exists an extension of $P$ to a sheaf $\oP$ on $\oJ\times\oJ$ such that the corresponding
Fourier-Mukai transform $\overline{\mathfrak F}:D^b(\oJ)\to D^b(\oJ)$ is an equivalence. After the first version
of this paper was completed, such an extension was constructed in \cite{cjacobians}. 
Then \eqref{eq:fouriermukai} is a composition of $\overline{\mathfrak F}$ and the direct image $j_*:D^b(J)\to D^b(\oJ)$ for the open
embedding $j:J\hookrightarrow\oJ$. Since $j_*$ is fully faithful, so is \eqref{eq:fouriermukai}. Thus our result is
natural assuming existence of $\overline{\mathfrak F}$.

(\theenumi)\stepcounter{enumi}
Compactified Jacobians appear as (singular)
fibers of the Hitchin fibration for the group $GL(n)$; therefore, our results can be interpreted as a kind of autoduality
of the Hitchin fibration. Conversely, some of our results can be derived from a theorem of E.~Frenkel and C.~Teleman
\cite{Teleman} (see Theorem \ref{th:hitchinpushforward}).
We explore this relation in more details in Section \ref{sc:Hitchin}.

(\theenumi)\stepcounter{enumi}
Recall that the curve $C$ is assumed to be integral with planar singularities. We assume integrality of $C$
to avoid working with stability conditions for sheaves on $C$. It is likely that our argument works without
this assumption if one fixes an ample line bundle on $C$ and defines the compactified Jacobian $\oJ$ to be the moduli
space of  semi-stable torsion-free sheaves of degree zero. Such generalization is natural in view of the previous remark,
because some fibers of the Hitchin fibration are compactified Jacobians of non-integral curves.

On the other hand, the assumption that $C$ has planar singularities is more important.
There are two reasons why
the assumption is natural. First of all, $\oJ$ is irreducible if and only if the singularities of $C$ are planar
(\cite{reducibility}); so if one drops this assumption, $J$ is no longer dense in $\oJ$. Secondly, only compactified
Jacobians of curves with planar singularities appear in the Hitchin fibration.

\end{remarks*}

\subsection*{Acknowledgments}

I would like to thank R.~Bezrukavnikov for stimulating my interest in this subject. This text was influenced by my
numerous
discussions with T.~Pantev, and I am very grateful to him for sharing his ideas.
I would also like to thank V.~Drinfeld, T.~Graber, C.~Teleman, J.~Starr, and J.~Wahl for their remarks and suggestions.

\section{Main results}
Fix a ground field $\kkk$. For convenience, let us assume that $\kkk$ is algebraically closed.
Let $C$ be an integral projective curve over $\kkk$. Denote by $J$ its Jacobian, that is,
$J$ is the moduli space of line bundles on $C$ of degree zero. Denote by $\oJ$ the compactified
Jacobian; in other words, $\oJ$ is the moduli space of torsion-free sheaves
on $C$ of generic rank one and degree zero.
(For a sheaf $F$ of generic rank one, the degree is $\deg(F)=\chi(F)-\chi(O_C)$.)

Let $P$ be the Poincar\'e bundle; it is a line bundle on $J\times\oJ$. Its fiber over $(L,F)\in (J\times\oJ)$
equals
\begin{equation}
P_{(L,F)}=\detrg(L\otimes F)\otimes\detrg (O_C)\otimes\detrg(L)^{-1}\otimes\detrg(F)^{-1}.
\end{equation}
More explicitly, we can write $L\simeq O(\sum a_ix_i)$ for a divisor $\sum a_ix_i$ supported by the smooth locus
of $C$, and then
$$P_{(L,F)}=\bigotimes(F_{x_i})^{\otimes a_i}.$$

From now on, we assume that $C$ has planar singularities; that is, the tangent space to $C$ at any point is at most two-dimensional.
Our main result is the computation of the direct image of $P$:

\begin{THEOREM}
$$Rp_{1,*}P=\det(H^1(C,O_C))\otimes O_\zeta[-g].$$ Here $O_\zeta$ is the structure sheaf of the neutral element $\zeta=[O_C]\in J$, and
$p_1:J\times\oJ\to J$ is the projection.
\label{th:pushforward}
\end{THEOREM}

Let us now view $P$ as a family of line bundles on $\oJ$ parametrized by $J$. For fixed $L\in J$, denote the corresponding
line bundle on $\oJ$ by $P_L$. In other words, $P_L$ is the restriction of $P$ to $\{L\}\times\oJ$. Applying base change, we can use Theorem \ref{th:pushforward}
to compute cohomology of $P_L$:

\begin{THEOREM}
\begin{enumerate}
\item If $L\not\simeq O_C$, then $H^i(\oJ,P_L)=0$ for any $i$;\label{th:cohomology1}

\item If $L=O_C$, then $P_L=O_\oJ$ and $H^i(\oJ,O_\oJ)=\bigwedge^i H^1(C,O_C)$. (The identification is described
more explicitly in Proposition \ref{pp:cohomology}.)\label{th:cohomology2}
\end{enumerate}
\qed
\label{th:cohomology}
\end{THEOREM}

Let $\Pic(\oJ)$ be the moduli space of line bundles on $\oJ$. The correspondence $L\mapsto P_L$ can be viewed
as a morphism $\rho:J\to\Pic(\oJ)$. Denote by $\Pic^0(\oJ)\subset\Pic(\oJ)$ the connected component of the identity
$[O_\oJ]\in\Pic(\oJ)$. In Section \ref{sc:autoduality}, we derive the following statement.

\begin{THEOREM}
$\rho$ gives an isomorphism $J\iso\Pic^0(\oJ)$.
\label{th:Jacobians}
\end{THEOREM}

\begin{remark*} Theorem \ref{th:Jacobians} answers the question raised in \cite{autoduality}.
Following \cite{GrothendieckVI}, set
\begin{equation}
\begin{aligned}
\Pic^\tau(\oJ)&=\{L\in\Pic(\oJ):L^{\otimes n}\in\Pic^0(\oJ)\text{ for some }n>0\},\\
\Pic^\sigma(\oJ)&=\{L\in\Pic(\oJ):L^{\otimes n}\in\Pic^0(\oJ)\text{ for some }n\text{ coprime to }\chr\kkk\}
\end{aligned}
\label{eq:tausigma}
\end{equation}
(if $\chr\kkk=0$, $\Pic^\sigma(\oJ)=\Pic^\tau(\oJ)$ by definition).
The main result of \cite{autoduality} is the Autoduality Theorem, which claims that
if all singularities of $C$ are double points, then $\rho:J\iso\Pic^0(\oJ)$ and $\Pic^0(\oJ)=\Pic^\tau(\oJ)$.
Theorem \ref{th:Jacobians} generalizes the first statement to curves with planar singularities; as for the second
statement, we show in Proposition \ref{pp:sigma} that $\Pic^0(\oJ)=\Pic^\sigma(\oJ)$. We do not know whether
$\Pic^\tau(\oJ)$ and $\Pic^\sigma(\oJ)$ coincide when $\chr(\kkk)> 0$ and $C$ has planar singularities.
\end{remark*}

Theorem~\ref{th:pushforward} can be reformulated in terms of the Fourier functor
$${\mathfrak F}:D^b(J)\to D^b(\oJ):\cF\mapsto Rp_{2,*}(p_1^*(\cF)\otimes P)$$
given by $P$. Recall that $D^b(J)$ stands for the (bounded) derived category of quasicoherent sheaves on $J$. 
The functor $\mathfrak F$ admits a left adjoint given by 
$${\mathfrak F}^\vee:D^b(\oJ)\to D^b(J):\cF\mapsto Rp_{1,*}(p_2^*(\cF)\otimes P^{-1})\otimes\det(H^1(C,O_C))^{-1}[g].$$
This formula relies on the computation of the dualizing sheaf on $\oJ$: see Corollary~\ref{co:dualizing}.

\begin{THEOREM}
\begin{enumerate}
\item
The composition ${\mathfrak F}^\vee\circ{\mathfrak F}$ is isomorphic to the identity functor.
\item ${\mathfrak F}$ is fully faithful. 
\end{enumerate}
\label{th:fourier}
\end{THEOREM}
\begin{proof}
The first statement follows from Theorem~\ref{th:pushforward} by base change. (This is completely analogous to the original argument 
of \cite[Theorem 2.2]{Mukai}.) This implies the second statement, because the functors ${\mathfrak F}^\vee$ and ${\mathfrak F}$ are adjoint.
\end{proof}

\begin{remark*} For simplicity, we considered a single curve $C$ in this section. However, all our results hold for
families of curves. Actually, we prove Theorem \ref{th:pushforward} for the universal family of curves
(Theorem \ref{th:main}); base change then implies that the statement holds for any family, and, in particular, for any single curve.
\end{remark*}

\section{Line bundles on a compactified Jacobian}

\begin{proposition}
Suppose $H^i(\oJ,P_L)\ne 0$ for some $i$. Then $(P_L)|_J\simeq O_J$.
\label{pp:naivevanishing}
\end{proposition}
\begin{proof}
Let $T\to J$ be the $\gm$-torsor corresponding to $(P_L)|_J$. One easily sees that $T$ is naturally an abelian group that is an extension of $J$ by $\gm$. The action of $J$ on $\oJ$ lifts to an action of $T$ on $P_L$, therefore, $T$ also
acts on $H^i(\oJ,P_L)$. Note that $\gm\subset T$ acts via the tautological character.

Let $V\subset H^i(\oJ,P_L)$ be an irreducible $T$-submodule. Since $T$ is commutative, $\dim(V)=1$. The action of $T$
on $V$ is given by a character $\chi:T\to\gm$. Since $\chi|_{\gm}= id$, we see that $\chi$ gives a splitting $T\simeq\gm\times J$. This implies the
statement.
\end{proof}

\begin{remark*} If $C$ is smooth, Proposition \ref{pp:naivevanishing} is equivalent to observation (vii) in
\cite[Section II.8]{Mumford}; however,
our proof uses a slightly different idea, which is better adapted to the singular case.
\end{remark*}

Let $C^0\subset C$ be the smooth locus of $C$.

\begin{corollary}
Suppose $H^i(\oJ,P_L)\ne 0$ for some $i$. Then $L|_{C^0}\simeq O_{C^0}$.
\label{co:naivevanishing2}
\end{corollary}
\begin{proof}
Fix a degree minus one line bundle $\ell$ on $C$. It defines an Abel-Jacobi map $$\alpha:C\to\oJ:c\mapsto\ell(c).$$
Here $\ell(c)$ can be defined as the sheaf of homomorphisms from the ideal sheaf of $c\in C$ to $\ell$.
Notice that $\alpha^*(P_L)\simeq L$ and $\alpha(C^0)\subset J$. Now Proposition
\ref{pp:naivevanishing} completes the proof.
\end{proof}

Set $$N=\{ L\in J:H^i(\oJ,P_L)\ne 0\text{ for some }i\}\subset J.$$ \
Clearly, $N\subset J$ is closed (by the Semicontinuity Theorem), and
$N=\supp(Rp_{1,*}P)$, where $p_1:J\times\oJ\to J$ is the projection (by base change).

\begin{corollary}
Let $g$ be the (arithmetic) genus of $C$ and $\tilde g$ be its geometric genus, that is, the genus of
its normalization. Then
$\dim(N)\le(g-\tilde g)$.
\label{co:naivevanishing}
\end{corollary}
\begin{proof}
Let $\nu:\tC\to C$ be the normalization, and let $\tJ$ be the Jacobian of $\tC$. The map $\nu^*:J\to\tJ$ is smooth and
surjective; its fibers have dimension $(g-\tilde g)$.

Denote by $\tN\subset\tJ$ the set of line bundles on $\tC$ that are trivial on $\nu^{-1}(C^0)\subset\tC$.
By Corollary \ref{co:naivevanishing2}, $\nu^*(N)\subset\tN$. Now it suffices to note that
$\tN$ is a countable set.
\end{proof}

\section{Moduli of curves}

Let $\cv=\cv_g$ be the moduli stack of integral projective curves $C$ of genus $g$ with planar singularities.
The following properties of $\cv$ are well known:

\begin{proposition}
$\cv$ is a smooth algebraic stack of finite type; $\dim(\cv)=3g-3$. \qed
\end{proposition}

\begin{remark*} Denote by $\cC$
the universal curve over $\cv$; that is, $\cC$ is the moduli stack of pairs $(C\in\cv,c\in C)$. One easily checks that $\cC$ is a smooth stack of dimension
$3g-2$. This is similar to the statement \eqref{th:jacobian2'} after Theorem \ref{th:jacobian}.
\end{remark*}

Consider the normalization $\tC$ of a curve $C\in\cv$, and let $\tg$ be the genus of $\tC$ (that is, the 
geometric genus of $C$). We need some results on the stratification of $\cv$ by geometric genus due to 
Tessier (\cite{Te}), Diaz and Harris (\cite{DH}), and Laumon (\cite{La}). Since our settings are somewhat different,
we provide the proofs.

Denote by $\cv^{(\tg)}\subset\cv$ the locus of curves $C\in\cv$ of geometric genus $\tg$. Note that we view $\cv^{(\tg)}$ simply as a subset of the set of points of $\cv$, rather than a substack.

\begin{proposition}
$\cv^{(\tg)}$ is a stratification of $\cv$:
$$\overline{(\cv^{(\tg)})}\subset\bigcup_{\gamma\le\tg}\cv^{(\gamma)}.$$
In particular, $\cv^{(\tg)}\subset\cv$ is locally closed.
\label{pp:stratification}
\end{proposition}
\begin{proof}
Let $\cS$ be the stack of birational morphisms $(\nu:\tC\to C)$, where $C\in\cv$, and $\tC$ is an integral projective curve of genus $\tg$
(with arbitrary singularity).
Consider the forgetful map
$$\pi:\cS\to\cv:(\nu:\tC\to C)\mapsto C.$$ Clearly,
$$\pi(\cS)\subset\bigcup_{\gamma\le\tg}\cv^{(\gamma)}.$$
Therefore, it suffices to show that $\pi$ is projective.

Let $\cS''$ be the stack of collections $(C,F,s)$, where $C\in\cv$, $F$ is a torsion-free sheaf on $C$ of generic rank
one and degree $g-\tg$, $s\in H^0(C,F)$.
Also, let $\cS'$ be the stack of collections $(C,F,s,\mu)$, where $(C,F,s)\in\cS''$ and
$\mu:F\otimes F\to F$ is such that $\mu(s\otimes s)=s$.
Consider
$$\cS\to\cS':(\nu:\tC\to C)\mapsto(C,\nu_*(O_\tC),1,\mu),$$
where $\mu$ is the product on the sheaf of algebras $\nu_*(O_\tC)$. This identifies $\cS$ and $\cS'$. The forgetful map
$$\cS'\to\cS'':(C,F,s,\mu)\mapsto(C,F,s)$$ is a closed embedding (essentially because $\mu$ is uniquely determined by
$\mu(s\otimes s)=s$).
Finally, the map
$$\cS''\to\cv:(C,F,s)\mapsto C$$
is projective.
\end{proof}

\begin{proposition}
$\codim(\cv^{(\tg)})\ge (g-\tg)$.
\label{pp:strata}
\end{proposition}
\begin{proof}
Let $\cS$ be as in the proof of Proposition \ref{pp:stratification}. Denote by $\cS^0$ the substack of morphisms
$(\nu:\tC\to C)\in\cS$ with smooth $\tC$; clearly, $\cv^{(\tg)}=\pi(\cS^0)$. Therefore, we need to show that
$\dim(\cS^0)\le 2g+\tg-3$.

Consider the morphism $$\tpi:\cS^0\to\cv_{\tg}:(\nu:\tC\to C)\mapsto\tC.$$
It suffices to show $\dim(\tpi^{-1}(\tC))\le 2(g-\tg)$ for any $\tC\in\cv_{\tg}$.
Fix $(\nu:\tC\to C)\in\cS^0$. Let us prove that the dimension of the tangent space
$T_\nu\tpi^{-1}(\tC)$ to $\tpi^{-1}(\tC)$ at this point is at most $2(g-\tg)$.

$T_\nu\tpi^{-1}(\tC)$ is isomorphic to the space of first-order deformations of $O_C$ viewed as a sheaf of subalgebras
of $\nu_* O_\tC$. This yields an isomorphism
$$
T_\nu\tpi^{-1}(\tC)=\{\text{differentiations }O_C\to\nu_*O_\tC/O_C\}=\Hom_{O_C}(\Omega_C,\nu_*O_\tC/O_C).$$
Now it suffices to notice that the fibers of the cotangent sheaf $\Omega_C$ are at most two-dimensional, and that the length of the sky-scraper sheaf
$\nu_*O_\tC/O_C$ equals $g-\tg$.
\end{proof}

\begin{remark*} By looking at nodal curves, one sees that $\codim(\cv^{(\tg)})=g-\tg$.
\end{remark*}

\section{Universal Jacobian}
Let
$\overline\cJ$ (resp. $\cJ\subset\overline\cJ$) be the relative compactified Jacobian (resp. relative Jacobian) of $\cC$ over $\cv$.
Here is the precise definition:

\begin{definition} For a scheme $S$, let $\hat\cJ_S$ be the following groupoid:
\begin{itemize}
\item Objects of $\hat\cJ_S$ are pairs $(C,F)$, where $C\to S$ is a flat family of integral projective curves with planar singularities
(that is, $C\in\cv_S$), and $F$ is a $S$-flat coherent sheaf on $C$ whose
restriction to the fibers of $C\to S$ is torsion free of generic rank one
and degree zero;
\item Morphisms $(C_1,F_1)\to(C_2,F_2)$ are collections $$(\phi:C_1\iso C_2,\ell,\Phi:F_1\iso\phi^*(F_2)\otimes_{O_S}\ell),$$ where $\phi$ is a morphism of $S$-schemes, and
$\ell$ is an invertible
sheaf on $S$.
\end{itemize}
As $S$ varies, groupoids $\hat\cJ_S$ form a pre-stack; let $\overline\cJ$ be the stack associated to it. Also, consider pairs $(C,F)$ where $C\in\cv_S$ and
$F$ is a line bundle on $C$ (of degree zero along the fibers of $S\to C$); such pairs form a sub-prestack of $\hat\cJ$; let $\cJ\subset\overline\cJ$ be
the associated stack.
\end{definition}

Clearly, $\cJ\subset\overline\cJ$ is an open substack.
The main properties of these stacks are summarized in the following theorem (\cite{irreducibility}):

\begin{theorem}[Altman, Iarrobino, Kleiman]
\begin{enumerate}
\item
\label{th:jacobian1}
$\overline p:\overline\cJ\to\cv$ is a projective morphism with irreducible fibers of dimension $g$;

\item
\label{th:jacobian2}
$\overline p$ is locally a complete intersection;

\item The restriction $p:\cJ\to\cv$ is smooth.
\end{enumerate}
\qed
\label{th:jacobian}
\end{theorem}

\begin{remark*}
By \cite[Corollary~B.2]{FGS}, \eqref{th:jacobian2} can be strengthened:
\begin{enumerate}
\renewcommand\theenumi{\roman{enumi}'}
\setcounter{enumi}{1}

\item \label{th:jacobian2'}
$\overline\cJ$ is smooth.
\renewcommand\theenumi{\roman{enumi}}
\end{enumerate}
Clearly, \eqref{th:jacobian2'} together with \eqref{th:jacobian1} imply \eqref{th:jacobian2}.
\end{remark*}

\begin{remark*}
The key step in the proof of \eqref{th:jacobian1} is Iarrobino's calculation (see \cite{Hilbert}):
\begin{equation}
\dim(\mathrm{Hilb}_k(\kkk[[x,y]]))=k-1,
\label{eq:hilbert}
\end{equation}
where $\mathrm{Hilb}_k(\kkk[[x,y]])$ is the Hilbert scheme of codimension $k$ ideals in $\kkk[[x,y]]$.
For other proofs of \eqref{eq:hilbert}, see \cite{Briancon}, \cite[Theorem 1.13]{Nakajima} and \cite{Baranovsky}. Also,
J.~Rego gives an alternative inductive proof of \eqref{th:jacobian1} in \cite{Rego}.
\end{remark*}

Denote by $\frj$ the rank $g$ vector bundle on $\cv$ whose fiber over $C\in\cv$ is $H^1(C,O_C)$. Alternatively, $\frj$
can be viewed as the bundle of (commutative) Lie algebras corresponding to the group scheme $p:\cJ\to\cv$. The
relative dualizing sheaf for $p$ then equals $\Omega^g_{\cJ/\cv}=p^*(\det(\frj)^{-1})$. It is easy to find the dualizing
sheaf for $\overline p:\overline\cJ\to\cv$:

\begin{corollary}
The relative dualizing sheaf $\omega_{\overline p}$ of $\overline p$ equals $\overline p^*(\det(\frj)^{-1})$.
\label{co:dualizing}
\end{corollary}
\begin{proof}
By Theorem \ref{th:jacobian}\eqref{th:jacobian2}, $\overline p$ is Gorenstein, so $\omega_{\overline p}$ is a line bundle.
Since $\omega_{\overline p}|_\cJ=\Omega^g_{\cJ/\cv}$, it suffices to check that $\codim(\overline\cJ-\cJ)\ge 2$.
But this is clear because a generic curve $C\in\cv$ is smooth (see Proposition \ref{pp:strata}).
\end{proof}

\section{Proof of Theorem \ref{th:pushforward}}

Consider the Poincar\'e bundle on $\cJ\times_\cv\overline\cJ$. We still denote it by $P$.

\begin{theorem} Let $p_1:\cJ\times_\cv\overline\cJ\to\cJ$ be the projection. Then
$$Rp_{1,*}P=(\Omega^g_{\cJ/\cv})^{-1}\otimes\zeta_*O_\cv[-g]=\zeta_*\det(\frj)[-g],$$
where $\zeta:\cv\to\cJ$ is the zero section.
\label{th:main}
\end{theorem}
\begin{proof}
Consider the dual $P^{-1}=\HOM(P,O)$ of $P$. (Actually $P^{-1}=(\nu\times id)^* P$,
where $\nu:\cJ\to\cJ$ is the involution $L\mapsto L^{-1}$.)
By  Corollary \ref{co:dualizing}, the dualizing sheaf of $p_1$ is isomorphic to $p_1^*\Omega^g_{\cJ/\cv}$.
Therefore,
\begin{equation}
R\HOM(Rp_{1,*}P,O_\cJ)=(Rp_{1,*}P^{-1})\otimes\Omega^g_{\cJ/\cv}[g]
\label{eq:duality}
\end{equation}
by Serre's duality.

Combining Corollary \ref{co:naivevanishing} and Proposition \ref{pp:strata}, we see that 
$$\codim(\supp(Rp_{1,*}P))\ge g.$$ 
By \eqref{eq:duality}, we see that both $Rp_{1,*}P$ and $R\HOM(Rp_{1,*}P,O_\cJ)[-g]$ are
concentrated in cohomological degrees from zero to $g$. It is now easy to see that $Rp_{1,*}P$
is concentrated in cohomological degree $g$, and that $R^gp_{1,*}P$
is a coherent Cohen-Macaulay sheaf of codimension $g$.

Next, notice that the restriction of $P$ to $\zeta(\cv)\times_\cv\overline\cJ$ is trivial. This provides
a map
\[\zeta^*(R^gp_{1,*}P)\to R^g\overline p_*(O_\cJ).\]
By Serre's duality, $R^g\overline p_* O_\cJ=\det(\frj)$. Now by adjunction, we obtain a morphism
\begin{equation}
R^gp_{1,*}P\to\zeta_*\det\frj.
\label{eq:adj}
\end{equation}
It remains to verify that \eqref{eq:adj} is an isomorphism.
Since \eqref{eq:adj} is an isomorphism over $\zeta(\cv)$ by construction, we need
to verify that $\supp(R^gp_{1,*}P)=\zeta(\cv)\subset\cJ$.
 
Let us check that $\supp(R^gp_{1,*}P)$ equals $\zeta(\cv)$ as a set.
As a set, $\supp(R^gp_{1,*}P)$ consists of pairs $(C,L)\in\cJ$ such that the line bundle $L$ on $C$ 
satisfies $H^g(\oJ,P_L)\ne 0$. In this case, $H^0(\oJ,P^{-1}_L)\ne 0$ by Serre's duality. Since $\oJ$
is irreducible, we see that the line bundle
$P^{-1}_L$ has a subsheaf isomorphic to $O_\oJ$. On the other hand, the line bundles
$P^{-1}_L=P_{L^{-1}}$ and $O_\oJ=P_{O}$ are algebraically equivalent, and therefore their
Hilbert polynomials coincide. Hence $P_L\simeq O_\oJ$. Finally, we can 
restrict $P_L$ to the image of the Abel-Jacobi map
(see the proof of Corollary~\ref{co:naivevanishing2}) to obtain $L\simeq O_C$.

To complete the proof, let us verify that $\supp(R^g_{1,*}P)=\zeta(\cv)$ as a scheme. Since $R^g_{1,*}P$
is Cohen-Macaulay of codimension $g$, it suffices to check the claim generically on $\zeta(\cv)$. We can thus
restrict ourselves to the open substack of smooth curves in $\cv$, and the claim reduces to \eqref{eq:pushforward}.
\end{proof}
\begin{remark*} The proof is similar to an argument of S.~Lysenko
(see proof of Theorem 4 in \cite{orthogonality}), see also D.~Mumford's proof
of the theorem in \cite[Section III.13]{Mumford}.
\end{remark*}

Using base change, one easily derives Theorem \ref{th:pushforward} from Theorem \ref{th:main}.

\section{Autoduality}\label{sc:autoduality}
Recall that the morphism $\rho:J\to\Pic_{\oJ}$ is given by $L\mapsto P_L$. Since
the tangent space to $J$ at $[O_C]$ (resp. to $\Pic(\oJ)$
at $[O_\oJ]$) equals $H^1(C,O_C)$ (resp. $H^1(\oJ,O_\oJ)$), the differential of $\rho$ at $[O_C]\in J$
becomes a linear operator
\[
d\rho:H^1(C,O_C)\to H^1(\oJ,O_\oJ).
\]
Let us give a more precise form of Theorem \ref{th:cohomology}\eqref{th:cohomology2}:
\begin{proposition}
$d\rho$ is an isomorphism, and the (super-commutative) cohomology algebra $H^\bullet(\oJ,O_\oJ)$ is freely
generated by $H^1(\oJ,O_\oJ)$.
\label{pp:cohomology}
\end{proposition}
\begin{proof}
Let $O_\zeta$ be the structure sheaf of the zero $[O_C]\in J$ viewed as a coherent sheaf on $J$ (it is a sky-scraper
sheaf of length one). Note that ${\mathfrak F}(O_\zeta)=O_\oJ$, where ${\mathfrak F}:D^b(J)\to D^b(\oJ)$ is the
Fourier transform of Theorem \ref{th:fourier}. Since $\mathfrak F$ is fully faithful, it induces an isomorphism
\begin{equation*}
\Ext^\bullet(O_\zeta,O_\zeta)\simeq\Ext^\bullet(O_\oJ,O_\oJ)=H^\bullet(\oJ,O_\oJ).
\end{equation*}
Finally, $J$ is smooth;
therefore, $\Ext^\bullet(O_\zeta,O_\zeta)=\bigwedge^\bullet H^1(C,O_C)$.
\end{proof}

Let us fix a line bundle $\ell$ on $C$ of degree minus one. It defines an Abel-Jacobi map $\alpha:C\to\oJ$, as in the proof
of Corollary~\ref{co:naivevanishing2}. We then obtain a morphism 
\[\alpha^*:\Pic(\oJ)\to\Pic(C):L\mapsto\alpha^* L.\] 
By construction, $\alpha^*$ is a left inverse of $\rho$ (cf. \cite[Proposition~2.2]{autoduality}).

\begin{remark*} Injectivity of $d\rho$ follows from the existence of the left inverse. Once injectivity is known, bijectivity
follows from the equality
$$\dim H^1(\oJ,O_\oJ)=\dim H^1(C,O_C)=g.$$
\end{remark*}

\begin{proof}[Proof of Theorem \ref{th:Jacobians}]
$\Pic(\oJ)$ is a group scheme of locally finite type 
(see \cite[Theorem~3.1]{GrothendieckV}, \cite[Theorem~9.4.8]{FAG}, or 
\cite[Corollary~(6.4)]{compactifying}). Set 
\begin{align*}
\Pic'(\oJ)&=(\alpha^*)^{-1}(J)=\{L\in\Pic(\oJ):\deg(\alpha^*L)=0\}\\
K&=\ker(\alpha^*)=\{L\in\Pic(\oJ):\alpha^*L\simeq O_C\}.
\end{align*}
Clearly, $K\subset\Pic'(\oJ)$ is closed, and $\Pic'(\oJ)\subset\Pic(\oJ)$ is both open and closed. 
The map
\[J\times K\to\Pic'(\oJ):(L_1,L_2)\mapsto\rho(L_1)\cdot L_2\]
is an isomorphism. Bijectivity of $d\rho$ (Proposition~\ref{pp:cohomology}) implies that $K$ is a disjoint union of points.
Therefore, the connected component of identity of $\Pic(\oJ)$ is contained in $\rho(J)$. 
Now it remains to notice that $J$ is connected.
\end{proof}

\begin{proposition}
$\Pic^\sigma(\oJ)=\Pic^0(\oJ)$ (where $\Pic^\sigma$ is defined in \eqref{eq:tausigma}).
\label{pp:sigma}
\end{proposition}
\begin{proof}
Consider $\overline p:\overline\cJ\to\cv$. It is a projective flat morphism with integral fibers (Theorem \ref{th:jacobian}); 
we can therefore construct the corresponding family of Picard schemes $\Pic(\overline\cJ/\cv)\to\cv$
(see the references in the proof
of Theorem~\ref{th:Jacobians}).
The family is separated and its fiber over $C\in\cv$ is $\Pic(\oJ_C)$.

Let us work in the smooth topology of $\cv$. Locally, we can choose a degree minus one line bundle
$\ell$ on the universal curve $\cC\to\cv$. As in the
proof of Theorem~\ref{th:Jacobians}, we then introduce a map
\[\alpha^*:\Pic(\overline\cJ/\cv)\to\Pic(\cC/\cv)\]
and substacks
$\Pic'(\overline\cJ/\cv)=(\alpha^*)^{-1}(\cJ)$ and $\cK=\ker(\alpha^*)$ such that
\[\Pic'(\overline\cJ/\cv)=\cJ\times_\cv\cK.\]

Let $\Pic^\sigma(\overline\cJ/\cv)\subset\Pic(\overline\cJ/\cv)$ be the substack whose
fiber over $C\in\cv$ is $\Pic^\sigma(\oJ_C)$. We have
\[\Pic^\sigma(\overline\cJ/\cv)=\cJ\times_\cv\cK^\sigma,\]
where 
\[\cK^\sigma=\{L\in\cK:L^{\otimes n}\simeq O\text{ for some $n$ coprime to $\chr\kkk$}\}.\]
By \cite[Theorem 2.5]{GrothendieckVI}, the map
\[\Pic(\overline\cJ/\cv)\to\Pic(\overline\cJ/\cv):L\mapsto L^{\otimes n}\]
is \'etale for all $n$ coprime to $\chr\kkk$. Therefore, $\cK^\sigma$ is \'etale over
$\cv$. 

Finally, the morphism $\cK^\sigma\to\cv$ is separated, and over the locus of smooth curves
$C\in\cv$, we have $\Pic^0(\oJ_C)=\Pic^\sigma(\oJ_C)$
by \cite[Corollary IV.19.2]{Mumford}. Therefore, $\cK^\sigma$ is the zero group scheme, and
$\Pic^\sigma(\overline\cJ/\cv)=\cJ$, as required.
\end{proof}

\section{Fibers of the Hitchin fibration}\label{sc:Hitchin}

Recall the construction of the Hitchin fibration \cite{Hitchin} (for $GL(n)$). Fix a smooth curve $X$
and an integer $n$.

\begin{definition}
A \emph{Higgs bundle} is a rank $n$ vector bundle $E$ on $X$ together with a \emph{Higgs field}
$A:E\to E\otimes\Omega_X$.
\end{definition}

Given a Higgs bundle $(E,A)$, consider the characteristic polynomial of $A$:
\begin{equation}
\det(\lambda I-A)=\lambda^n+a_1\lambda^{n-1}+\dots+a_n; \qquad a_i\in H^0(X,\Omega^{\otimes i}_X).
\label{eq:characteristic}
\end{equation}
The zero locus of 
\eqref{eq:characteristic} is a curve $C\subset T^*X$: the \emph{spectral curve} of $A$. Higgs bundle $(E,A)$ gives rise
to a coherent sheaf $F$ on $C$; informally, $F$ is the
`sheaf of co-eigenspaces': its fiber over a point $(x,\mu)\in T^*X$ is the co-eigenspace 
$$\mathop{\mathrm{coker}}(A(x)-\mu:E_x\to E_x\otimes\Omega_{X,x}).$$ Here $x\in X$, $\mu\in\Omega_{X,x}$.

\begin{proposition}
\label{pp:hitchin}
\begin{enumerate}
\item $F$ is a torsion-free sheaf on $C$ whose fiber at a generic point of $C$ has length equal to
the multiplicity of the corresponding component of $C$. In particular, if $C$ is reduced, $F$ is a 
torsion-free sheaf of generic rank one. \label{pp:hitchin1}

\item Fix a spectral curve $C$ (that is, fix a polynomial \eqref{eq:characteristic}). Then $(E,A)\mapsto F$
is a one-to-one correspondence between Higgs bundles with spectral curve $C$ and sheaves $F$ as in \eqref{pp:hitchin1}.
\end{enumerate}
\qed
\end{proposition}

Given $F$, $E$ is reconstructed as the push-forward of $F$ with respect to $C\to X$. Therefore, $F$ and $E$ 
have equal Euler characteristics. We have $\chi(O_C)=n^2\chi(O_X)=n^2(1-g)$, where $g$ is the genus of $X$. Hence
$\deg(F)=0$ if and only if $\deg(E)=n(n-1)(1-g)$, where 
$g$ is the genus of $X$. (Recall that 
$\deg(F)=\chi(F)-\chi(O_C)$.) Also, note that $(E,A)$ is (semi)stable if
and only if $F$ is (semi)stable. 
If $C$ is integral, $F$ has generic rank one and stability is automatic.

Let $\Higgs$ be the moduli space of semi-stable Higgs bundles $(E,A)$ with $\rk(E)=n$ and $\deg(E)=n(n-1)(1-g)$.
Also, let $\SCurves$ be the space of spectral curves 
$C\subset T^*X$; explicitly, $\SCurves$ is the space of coefficients $(a_1,\dots,a_n)$ of \eqref{eq:characteristic}:
$$\SCurves=\prod_{i=1}^n H^0(X,\Omega^{\otimes i}_X).$$ Finally, let $\SCurves'\subset\SCurves$ be the locus of 
integral spectral curves $C\subset T^*X$.

The correspondence $(E,A)\mapsto C$ gives a map $h:\Higgs\to\SCurves$ (the \emph{Hitchin fibration}). For $C\in\SCurves$,
the fiber $h^{-1}(C)$ is the space of Higgs bundles with spectral curve $C$; Proposition \ref{pp:hitchin} 
identifies $h^{-1}(C)$ with the moduli space of semi-stable coherent sheaves $F$ on $C$ that satisfy 
Proposition \ref{pp:hitchin}\eqref{pp:hitchin1} and have degree zero. In other words, the fiber is the compactified 
Jacobian of $C$. 

The results of this paper can be applied to integral spectral curves $C\in\SCurves'$. For instance, Theorem 
\ref{th:cohomology}\eqref{th:cohomology2} implies that $$H^i(h^{-1}(C),O)=\bigwedge\nolimits^i H^1(C,O_C).$$ Actually, applying
the relative version of Theorem \ref{th:cohomology}\eqref{th:cohomology2} to the universal family of spectral curves,
we obtain an isomorphism
\begin{equation}
(R^ih_* O_{\Higgs})|_{\SCurves'}=\Omega^i_{\SCurves'},
\label{eq:hitchinpushforward}
\end{equation}
where we used the symplectic form on $T^*X$ to identify $H^1(C,O_C)$ with the cotangent space to 
$C\in\SCurves'$. Recently, E.~Frenkel and C.~Teleman proved that 
the isomorphism \eqref{eq:hitchinpushforward} can be extended to the space of all spectral curves:

\begin{theorem} There is an isomorphism
$$
R^ih_* O_{\Higgs}=\Omega^i_{\SCurves}.
$$
\label{th:hitchinpushforward}
\qed
\end{theorem}

When $i=0,1$, Theorem \ref{th:hitchinpushforward} is proved by N.~Hitchin (\cite[Theorems 6.2 and 6.5]{Hitchin});
the general case is announced in \cite{Teleman}.

\begin{remarks*}
\setcounter{enumi}{1}
(\theenumi)\stepcounter{enumi}
In \cite{Hitchin}, N.~Hitchin works with the Hitchin fibration for the group $SL(2)$, but his argument
can be used to compute $R^ih_* O_{\Higgs}$ for arbitrary $n$ (still assuming $i=0,1$). Actually, essentially the same
argument computes $R^i{\overline p}_* O_{\overline\cJ}$ for $i=0,1$. (Recall that $\overline p:\overline\cJ\to\cv$ is the 
universal compactified Jacobian over the moduli stack of curves $\cv$.)

(\theenumi)\stepcounter{enumi}
In \cite{Teleman}, Theorem \ref{th:hitchinpushforward} is stated for the Hitchin fibration for arbitrary group, not just
$GL(n)$. 

(\theenumi)\stepcounter{enumi} One can derive some of our results from Theorem \ref{th:hitchinpushforward}, at least
for integral curves $C$ that appear as spectral curves of the Hitchin fibration. Indeed, for such $C\in\SCurves'$,
Theorem \ref{th:hitchinpushforward} implies Theorem \ref{th:cohomology}\eqref{th:cohomology2}. In turn, this implies
Theorem \ref{th:jacobian}. Also, one can easily derive from Theorem \ref{th:cohomology}\eqref{th:cohomology2} that the
isomorphism of Theorem \ref{th:pushforward} exists on some neighborhood $U$ of $\zeta\in J$, so 
Theorem \ref{th:cohomology}\eqref{th:cohomology1} holds for $L\in U$. Similarly, we see that $P$ defines a fully faithful
Fourier-Mukai transform from $D^b(U)$ to $D^b(\oJ)$. 
\end{remarks*}

\bibliographystyle{abbrv}
\bibliography{jacobians}
\end{document}